# SECOND CLASS PARTICLES AND CUBE ROOT ASYMPTOTICS FOR HAMMERSLEY'S PROCESS

BY ERIC CATOR AND PIET GROENEBOOM

*Delft University of Technology and Vrije Universiteit Amsterdam*

We show that, for a stationary version of Hammersley's process, with Poisson sources on the positive $x$-axis and Poisson sinks on the positive $y$-axis, the variance of the length of a longest weakly North–East path $L(t,t)$ from $(0,0)$ to $(t,t)$ is equal to $2\mathbb{E}(t - X(t))_+$, where $X(t)$ is the location of a second class particle at time $t$. This implies that both $\mathbb{E}(t - X(t))_+$ and the variance of $L(t,t)$ are of order $t^{2/3}$. Proofs are based on the relation between the flux and the path of a second class particle, continuing the approach of Cator and Groeneboom [*Ann. Probab.* **33** (2005) 879–903].

**1. Introduction.** In an influential paper Kim and Pollard [8] show that in many statistical contexts we are confronted with estimators which converge at rate $n^{1/3}$ instead of the usual rate $n^{1/2}$ and that in this situation the limit distribution is nonnormal. They call this phenomenon "cube root asymptotics." A prototype of such an estimator is the maximum likelihood estimator of a decreasing density, which converges locally at rate $n^{1/3}$ after rescaling to the (almost surely unique) location of the maximum of Brownian motion minus a parabola. The characterization of this limit distribution in terms of Airy functions was given in [6].

It has been conjectured that the asymptotics for longest increasing subsequences, which can be analyzed by studying longest North–East paths of Hammersley's process, is related to these cube root phenomena in estimation theory and, in particular, that it should be possible to derive the asymptotics along similar lines. However, up till now, the cube root limit theory for longest increasing subsequences and longest North–East paths has been based on certain analytic relations, involving Toeplitz determinants; see, for example, [2] and [3].

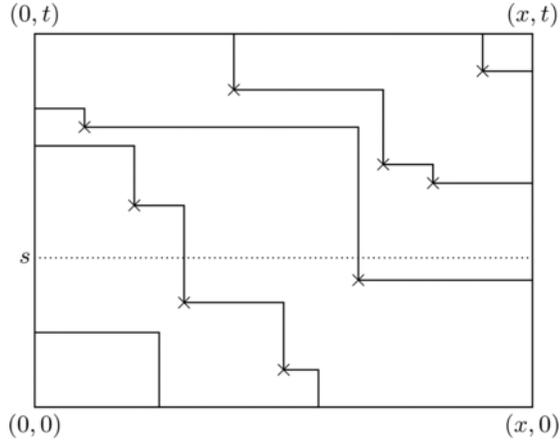

Fig. 1. *Space–time paths of the Hammersley's process, with sources and sinks.*

In this paper we will work with Hammersley's process with sources and sinks, as defined in [4]. We will give a short description here, based on Figure 1. We consider the space–time paths of particles that started on the $x$-axis as *sources*, distributed according to a Poisson distribution with parameter $\lambda$, and we consider the $t$-axis as a time axis. In the positive quadrant we have a Poisson process of what we call $\alpha$-points (denoted in Figure 1 by $\times$), which will have intensity 1, unless otherwise specified. On the $t$-axis (which also sometimes will be called $y$-axis) we have a Poisson process of what we call *sinks* of intensity $1/\lambda$. The three Poisson processes are independent.

At the time an $\alpha$-point appears, the particle immediately to the right of it jumps to the location of the $\alpha$-point. At the time a sink appears, the leftmost particle disappears. To know the particle configuration at time $s$, we intersect a line at time $s$ with the space–time paths. The counting process of the particle configuration at time $t$ is denoted by $L_\lambda(\cdot, t)$, where we start counting at the first sink on the $t$-axis up to $(0, t)$, and continue counting on the halfline $(0, \infty) \times \{t\}$, so $L_\lambda(x, t)$ equals the total number of sinks in the segment $\{0\} \times [0, t]$ plus the number of crossings of space–time paths of the segment $[0, x] \times \{t\}$.

The total number of space–time paths in $[0, x] \times [0, t]$ is called the *flux* at $(x, t)$. It is in fact equal to $L_\lambda(x, t)$. If $\lambda = 1$, we will denote $L_1(x, t)$ just by $L(x, t)$, unless this can cause confusion. The flux $L_\lambda(x, t)$ equals the length of a longest *weakly* NE (North–East) path from $(0, 0)$ to $(x, t)$, where "weakly" means that we are allowed to pick up either sources from the $x$-axis or sinks from the $t$-axis, before we start picking up $\alpha$-points. To see this from Figure 1, trace back a longest weakly NE path from $(x, t)$ to $(0, 0)$, and note that one will pick up exactly one $\alpha$-point or one source or one sink from each

space–time path. Note that, if $0 < x < y$, $L_\lambda(y,t) - L_\lambda(x,t)$ is the number of particles (or crossings of space–time paths) on the segment $[x,y] \times \{t\}$.

A heuristic argument for the cube root behavior of the fluctuation of the length of a longest weakly NE path for the stationary Hammersley process runs as follows. Suppose, for simplicity, that $\lambda = 1$. A longest weakly NE path with length $L(t,t) = L_1(t,t)$ can pick up points from either $x$- or $y$-axis before starting on a strictly NE path to $(t,t)$. Furthermore, let, for $-t \leq z \leq t$,

(1.1) $\quad N(z) = \begin{cases} \text{number of sources in } [0,z] \times \{0\}, & \text{if } z \geq 0, \\ \text{number of sinks in } \{0\} \times [0,|z|], & \text{if } z \leq 0, \end{cases}$

and

(1.2) $A_t(z) = \begin{cases} \text{length of longest strictly NE path from } (z,0) \text{ to } (t,t), \\ \quad \text{if } z \geq 0, \\ \text{length of longest strictly NE path from } (0,|z|) \text{ to } (t,t), \\ \quad \text{if } z < 0. \end{cases}$

Note that the processes $A_t$ and $N$ are independent and that

(1.3) $\qquad\qquad L(t,t) = \sup\{N(z) + A_t(z) : -t \leq z \leq t\}.$

The process $z \mapsto t^{-1/3}\{N(zt^{2/3}) - |z|t^{2/3}\}$, $|z| \leq t^{1/3}$, converges in the topology of uniform convergence on compacta to two-sided Brownian motion, originating from zero. As will be shown below, the expectation of $t^{-1/3}A_t(zt^{2/3})$ has an asymptotic upper bound (as $t \to \infty$) of the form

$$2t^{2/3} - |z|t^{1/3} - \tfrac{1}{4}z^2,$$

which is seen by taking expectations and optimizing the choice of $\lambda$ in the inequality in Lemma 4.1. This suggests that the distance to zero of the exit point where the longest path leaves either $x$- or $y$-axis cannot be of larger order than $t^{2/3}$, since otherwise the Brownian motion cannot cope with the downward parabolic drift $-\tfrac{1}{4}z^2$, temporarily assuming that the fluctuation of $A_t(z)$ is of order $O_p(t^{1/3})$. The latter fact we know to be true from the analytic approach, not used in our probabilistic approach. On the other hand, we will derive this in Section 7; see (7.7).

The limit behavior of the exit point can be compared to the behavior of the location of the maximum of Brownian motion minus a parabola, which plays a key role in the asymptotics for the cube root estimation theory, mentioned above. The crucial difference, however, is that the exit point is the location of the maximum of the sum of two independent processes instead of the maximum of just one process.

We note here that the $n^{1/3}$ convergence in estimation theory (so *slower* convergence than the usual $n^{1/2}$-convergence) corresponds to the $t^{2/3}$ order of the distance to zero of the exit point, which (after a time-reversal argument, based on Burke's theorem for Hammersley's process) can be called "super-diffusive" behavior of a second class particle. The slower convergence in estimation theory is caused by the fact that the estimators have an interpretation in terms of the *location* of a maximum, just as the exit point for a longest weakly NE path in Hammersley's process.

The key relation which allows us to make the heuristic argument above rigorous is Theorem 2.1 of Section 2, which, in combination with (time-reversal) results of Section 3, tells us that $\text{Var}(L(t,t)) = 2\mathbb{E}Z(t)_+$, where $Z(t)$ is the rightmost point where a longest weakly NE path leaves the $x$-axis; see (3.2), Section 3. It is shown in Section 4 that this implies $\mathbb{E}Z(t)_+ = O(t^{2/3})$. In Section 5 we compare longest strictly NE paths with longest weakly NE paths and obtain a bound on the difference between the lengths of these paths. This allows us to also obtain a lower bound for $\mathbb{E}Z(t)_+$ in Section 6. Finally, in Section 7 we discuss tightness results for the original Hammersley process without sources and sinks, in connection with results of Seppäläinen [9].

Our methods heavily rely on the ideas developed in [4], which concern, in particular, the difference in behavior below and above the path of a second class particle and Burke's theorem for Hammersley's process, which enables us to use time-reversal and reflection.

**2. Variance of the flux and location of a second class particle.** We will need the concept's second class particle and dual second class particle, which also play an important role in later sections. A "normal" second class particle is created by putting an extra source at $(0,0)$ (thus effectively removing the first sink), and a dual second class particle is created by putting an extra sink at $(0,0)$, thus effectively removing the first source. Define $X(t)$ as the location at time $t$ of a second class particle in a stationary Hammersley process with source and sink intensity equal to 1 (the symmetric case), and $X'(t)$ as the location at time $t$ of a dual second class particle for this case.

As explained in [4], a "normal" second class particle $X(t)$ jumps to the previous position of the ordinary ("first class") particle that exits through the first sink at the time of exit, and successively jumps to the previous positions of particles directly to the right of it, at times where these particles jump to a position to the left of the second class particle. The concept of a dual second class particle was also considered in [4], but there it is seen as a second class particle for the process "moving from left to right." Figure 2 shows the trajectories of a second class particle and a dual second class particle. Note that we always have $X(t) \leq X'(t)$, which is evident from Figure 2.

Now consider a stationary Hammersley process with $\alpha$-intensity 1, source-intensity $\lambda$ and sink-intensity $1/\lambda$. Fix $x, t > 0$ and consider the flux $L_\lambda(x,t)$.

Denote $X_\lambda(t)$ and $X'_\lambda(t)$ as the locations at time $t$ of a second class particle and a dual second class particle, respectively. We use the subscript $\lambda$ to indicate that the distribution of the location of the (dual) second class particle depends on $\lambda$. If $\lambda = 1$, the subscript is suppressed. We have the following result:

THEOREM 2.1.
$$\operatorname{Var}(L_\lambda(x,t)) = -\lambda x + \frac{t}{\lambda} + 2\lambda \mathbb{E}(x - X_\lambda(t))_+.$$

REMARK 2.1. A similar relation between the variance of the flux and the location of a second class particle has been proved for totally asymmetric simple exclusion processes (TASEP) in [5].

REMARK 2.2. Note that taking $\lambda = \sqrt{t/x}$ yields
$$\operatorname{Var}(L_\lambda(x,t)) = 2\lambda \mathbb{E}(x - X_\lambda(t))_+.$$

PROOF OF THEOREM 2.1. For notational clarity we use the four wind directions $N, E, S$ and $W$ to denote the number of crossings of the four respective sides of the rectangle $[0, x] \times [0, t]$ [so $L_\lambda(x,t) = N + W$]. Clearly, $S + E = N + W$. We also know from Burke's theorem for Hammersley's process (see [4]) that $N$ and $E$ are independent, just like $S$ and $W$. This means that
$$\operatorname{Var}(L_\lambda(x,t)) = \operatorname{Var}(W + N)$$

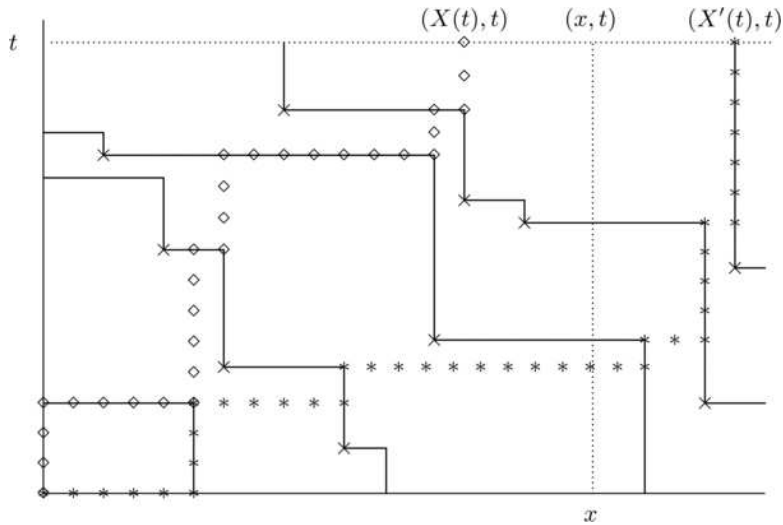

FIG. 2. *Trajectories of $(X(t), t)$ and $(X'(t), t)$.*

$$\begin{aligned}
&= \operatorname{Var}(W) + \operatorname{Var}(N) + 2\operatorname{Cov}(W, N) \\
&= \operatorname{Var}(W) + \operatorname{Var}(N) + 2\operatorname{Cov}(S + E - N, N) \\
&= \operatorname{Var}(W) - \operatorname{Var}(N) + 2\operatorname{Cov}(S, N) \\
&= \frac{t}{\lambda} - \lambda x + 2\operatorname{Cov}(S, N).
\end{aligned} \tag{2.1}$$

We want to investigate $\operatorname{Cov}(S, N)$. It turns out that we can do this by varying the source-intensity appropriately. For $\varepsilon > 0$, we define a source-intensity of $\lambda + \varepsilon$. The sinks remain a Poisson process with intensity $1/\lambda$. We denote expectations with respect to this new source intensity by $\mathbb{E}_\varepsilon$. Define

$$a_n = \mathbb{E}_\varepsilon(N | S = n).$$

Note that $a_n$ does, in fact, not depend on $\varepsilon$, since we condition on the number of sources in $[0, x]$, and the sources outside this interval do not influence $N$. Then

$$\begin{aligned}
\left.\frac{\partial}{\partial \varepsilon}\right|_{\varepsilon=0} \mathbb{E}_\varepsilon(N) &= \left.\frac{\partial}{\partial \varepsilon}\right|_{\varepsilon=0} \sum_{n=0}^\infty \frac{(x(\lambda+\varepsilon))^n}{n!} e^{-x(\lambda+\varepsilon)} a_n \\
&= \frac{1}{\lambda} \sum_{n=0}^\infty \frac{(x\lambda)^n}{n!} e^{-x\lambda} a_n \cdot n - x \sum_{n=0}^\infty \frac{(x\lambda)^n}{n!} e^{-x\lambda} a_n \\
&= \frac{1}{\lambda} \mathbb{E}(NS) - x \mathbb{E}(N).
\end{aligned}$$

This shows that

$$\begin{aligned}
\operatorname{Cov}(N, S) &= \mathbb{E}(NS) - \mathbb{E}(N)\mathbb{E}(S) \\
&= \lambda \left.\frac{\partial}{\partial \varepsilon}\right|_{\varepsilon=0} \mathbb{E}_\varepsilon(N),
\end{aligned} \tag{2.2}$$

where we use that $\mathbb{E}(S) = \lambda x$.

We will calculate this derivative in the following manner. Fix, independently, a Poisson process of intensity 1 of $\alpha$-points in $(0, \infty)^2$, a Poisson process of sources of intensity $\lambda$ on the $x$-axis and a process of sinks of intensity $1/\lambda$ on the $t$-axis. Now we add an independent Poisson process of intensity $\varepsilon$ to the process of sources. Define $N_\varepsilon$ as the number of crossings of the North-side [i.e., $(0, x) \times \{t\}$] for the process with the added sources.

Note that if we add an extra source at $(z, 0)$, then $N$ increases by 1 if and only if $X_\lambda(t; z) < x$, where $X_\lambda(t; z)$ is the location of a second class particle at time $t$, which started at $(z, 0)$. We denote $X_\lambda(t) = X_\lambda(t; 0)$. This means that

$$\mathbb{E}(N_\varepsilon) = \mathbb{E}(N_0) + \varepsilon \int_0^x \mathbb{E}(\mathbb{1}_{\{X_\lambda(t;z) < x\}}) \, dz + O(\varepsilon^2). \tag{2.3}$$

Therefore, by using the stationarity of the Hammersley process,

$$\operatorname{Cov}(N, S) = \lambda \frac{\partial}{\partial \varepsilon}\bigg|_{\varepsilon=0} \mathbb{E}(N_\varepsilon)$$

$$= \lambda \int_0^x \mathbb{E}(\mathbb{1}_{\{X_\lambda(t;z) < x\}})\, dz$$

(2.4)
$$= \lambda \int_0^x \mathbb{P}(X_\lambda(t) < x - z)\, dz$$

$$= \lambda \int_0^x \mathbb{P}(x - X_\lambda(t) > z)\, dz$$

$$= \lambda \mathbb{E}(x - X_\lambda(t))_+.$$

Combining this with (2.1) gives

$$\operatorname{Var}(L_\lambda(x,t)) = -\lambda x + \frac{t}{\lambda} + 2\lambda \mathbb{E}(x - X_\lambda(t))_+. \qquad \square$$

Now consider a stationary Hammersley process with source (and sink) intensity 1. We denoted the flux of this process at $(x,t)$ by $L(x,t)$, and the location of a (dual) second class particle at time $t$, which started at $(0,0)$, by $X(t)$ [resp. $X'(t)$]. Note that under the map

$$(x,t) \mapsto (x/\lambda, \lambda t),$$

a stationary process with source intensity 1 gets transformed into a stationary process with source intensity $\lambda$ (and corresponding sink intensity $1/\lambda$). This rescaling argument shows

(2.5) $$\lambda X_\lambda(t) \stackrel{\mathcal{D}}{=} X(t/\lambda) \quad \text{and} \quad \lambda X'_\lambda(t) \stackrel{\mathcal{D}}{=} X'(t/\lambda),$$

where $\stackrel{\mathcal{D}}{=}$ denotes equality in distribution.

We would like to bound $\operatorname{Var}(L_\lambda(t,t))$ in terms of $\operatorname{Var}(L(t,t))$ in the case where $\lambda \geq 1$. Using Theorem 2.1 and (2.5), we get, using the inequality $(A+B)_+ \leq A_+ + B_+$,

$$\operatorname{Var}(L_\lambda(t,t)) = -\lambda t + t/\lambda + 2\mathbb{E}(\lambda t - t/\lambda + t/\lambda - X(t/\lambda))_+$$

$$\leq (\lambda - 1/\lambda)t + 2\mathbb{E}(t/\lambda - X(t/\lambda))_+$$

$$= (\lambda - 1/\lambda)t + \operatorname{Var}(L(t/\lambda, t/\lambda)).$$

If we show the intuitively clear result that, for $\lambda \geq 1$,

(2.6) $$\operatorname{Var}(L(t/\lambda, t/\lambda)) \leq \operatorname{Var}(L(t,t)),$$

we have proved that

(2.7) $$\operatorname{Var}(L_\lambda(t,t)) \leq (\lambda - 1/\lambda)t + \operatorname{Var}(L(t,t)).$$

We can show (2.6) by noting that Theorem 2.1 for $\lambda = 1$ is equivalent to

$$\text{Var}(L(x,t)) = -x + t + 2 \int_0^x \mathbb{P}(X(t) \leq z) \, dz. \tag{2.8}$$

Define

$$\phi(x,t) = \text{Var}(L(x,t)).$$

Clearly, $\phi$ is symmetric, since the source and sink intensities are equal, which gives reflection symmetry of the process.

Furthermore, (2.8) shows that $\phi$ is a continuously differentiable function, with

$$\partial_1 \phi(x,t) = -1 + 2\mathbb{P}(X(t) \leq x).$$

If we can show that $\mathbb{P}(X(t) \leq t) \geq 1/2$, we would have proved (2.6), since $\partial_1 \phi(t,t) \geq 0$ for a symmetric function $\phi$ implies that $\phi(t,t)$ is increasing in $t$.

Since reflecting Hammersley's process in the diagonal preserves the distribution, while interchanging the trajectories of $X'$ and $X$, we know that

$$\mathbb{P}(X(t) > x) = \mathbb{P}(X'(x) < t) \leq \mathbb{P}(X(x) < t).$$

Choosing $t = x$, we see that

$$\mathbb{P}(X(t) > t) \leq \mathbb{P}(X(t) < t) = \mathbb{P}(X(t) \leq t),$$

which shows that $\mathbb{P}(X(t) \leq t) \geq 1/2$. As noted before, this proves (2.6).

**3. Connection between second class particles and exit points.** As has already been noted in the Introduction, we can view the flux $L_\lambda(x,t)$ in two ways: it is the number of space–time paths in the square $[0,x] \times [0,t]$, but is also the length of the longest weakly NE path from $(0,0)$ to $(x,t)$, where "weakly NE" means that we are allowed to pick up sources or sinks, as well as $\alpha$-points, as long as we are going North–East. To work with this latter representation, which we will mainly use in the symmetric case when both the source- and the sink-intensity are 1, we define, for $-t \leq z \leq t$, $N(z)$ and $A_t(z)$ by (1.1) and (1.2). Remember that the processes $A_t$ and $N$ are independent and that

$$L(t,t) = \sup\{N(z) + A_t(z) : -t \leq z \leq t\}. \tag{3.1}$$

Another important aspect of this representation is the location at which a longest path leaves either the $x$-axis or the $y$-axis. Define

$$Z(t) = \sup\{z \in [-t,t] : N(z) + A_t(z) = L(t,t)\} \tag{3.2}$$

and

$$Z'(t) = \inf\{z \in [-t,t] : N(z) + A_t(z) = L(t,t)\}. \tag{3.3}$$

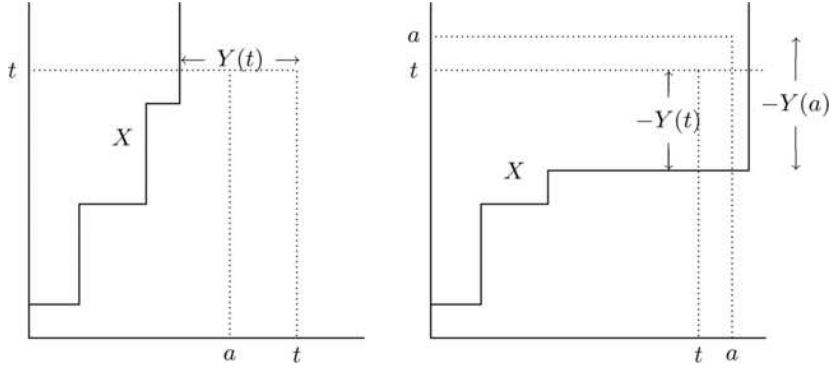

Fig. 3. *Relating $X$ and $Y$.*

We call $Z(t)$ and $Z'(t)$ *exit points* for a longest path, since there exist longest paths that leave the axis on $(Z(t), 0)$ [or $(0, -Z(t))$] or on $(Z'(t), 0)$ [or $(0, -Z'(t))$]. From this definition and using the symmetry of the situation, we can see that

(3.4) $$Z'(t) \leq Z(t) \quad \text{and} \quad Z(t) \stackrel{\mathcal{D}}{=} -Z'(t).$$

We will need another link between the two representations. We have defined $X(t)$ and $X'(t)$ as the position at time $t$ of a second class particle, respectively dual second class particle, that starts at $(0,0)$. Now define

$$Y(t) = \begin{cases} t - X(t), & \text{if } X(t) \leq t, \\ \inf\{s \geq 0 : X(s) \geq t\} - t, & \text{if } X(t) > t, \end{cases}$$

and

$$Y'(t) = \begin{cases} t - X'(t), & \text{if } X'(t) \leq t, \\ \inf\{s \geq 0 : X'(s) \geq t\} - t, & \text{if } X'(t) > t. \end{cases}$$

Since $X'(t) \geq X(t)$, we have $Y'(t) \leq Y(t)$. Figure 3 shows the relation between $X$ and $Y$.

It also shows two relations which we will be important later:

(3.5)
$$\begin{aligned}
a < t \quad &\Longrightarrow \quad \{X(t) < a\} = \{Y(t) > t - a\} \quad \text{and} \\
&\qquad \{X'(t) < a\} = \{Y'(t) > t - a\}, \\
a > t \quad &\Longrightarrow \quad \{X(t) > a\} = \{Y(a) < t - a\} \quad \text{and} \\
&\qquad \{X'(t) > a\} = \{Y'(a) < t - a\}.
\end{aligned}$$

Now consider Figure 4. The left picture shows a realization of the Hammersley process and two longest weakly NE paths, corresponding to $Z(t)$ and $Z'(t)$. The right picture shows the same realization, but now reflected in the point $(\frac{1}{2}t, \frac{1}{2}t)$.

Note that the longest paths become trajectories of a second class and a dual second class particle in the reflected process, and that $Z(t)$ corresponds to $Y(t)$, while $Z'(t)$ corresponds to $Y'(t)$. Burke's theorem in [4] states that the reflected process is also a realization of the stationary Hammersley process, so that we can indeed conclude that

$$Z(t) \stackrel{\mathcal{D}}{=} Y(t) \quad \text{and} \quad Z'(t) \stackrel{\mathcal{D}}{=} Y'(t).$$

In particular, this means, using Theorem 2.1 and noting that $(t - X(t))_+ = Y(t)_+$, that

(3.6) $$\operatorname{Var}(L(t,t)) = 2\mathbb{E}Z(t)_+.$$

**4. $EZ(t)_+$ is of order $O(t^{2/3})$.** We wish to control the exit point $Z(t)$. We will do this by considering an auxiliary Hammersley process $L_\lambda$, coupled to the original one, by thickening the sources to a Poisson process of intensity $\lambda \geq 1$ and thinning the sinks to a Poisson process of intensity $1/\lambda$. The process $A_t$ then satisfies the following inequality.

LEMMA 4.1. *Let $\lambda \geq 1$ and define $L_\lambda(x,t)$ as the flux of $L_\lambda$ at $(x,t)$. Then, for $0 \leq z \leq t$,*

$$A_t(z) \leq L_\lambda(t,t) - L_\lambda(z,0).$$

PROOF. It is clear that a strictly NE path from $(z, 0)$ to $(t, t)$ is shorter than a longest weakly NE path from $(z, 0)$ to $(t, t)$, where this path is allowed to either pick up sources of $L_\lambda$ on $[z, t] \times \{0\}$, or pick up crossings of $L_\lambda$ with $\{z\} \times [0, t]$. However, this longest weakly NE path is equal to the number of space–time paths in $[z, t] \times [0, t]$ of $L_\lambda$, which, in turn, is equal to $L_\lambda(t, t)$ minus the number of sources on $[0, z] \times \{0\}$. $\square$

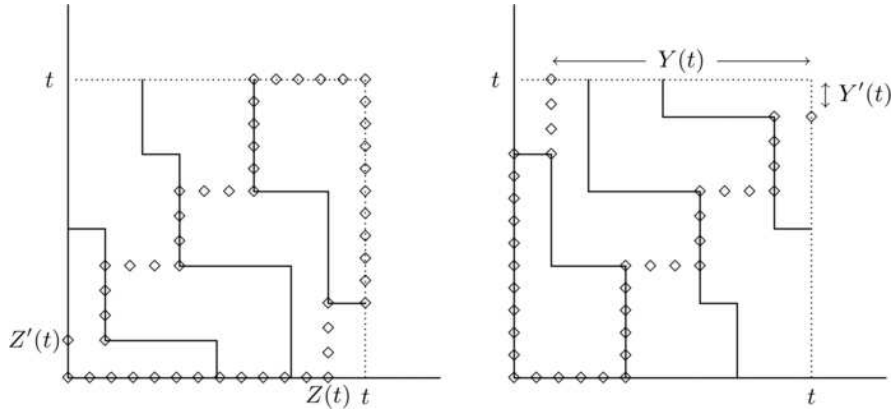

FIG. 4. *Longest path is distributed as trajectory of a second class particle.*

We can now show the following theorem. We use the notation $a(x) \lesssim b(x)$ if there exists a constant $M$ such that, for all parameters $x$, $a(x) \leq Mb(x)$.

THEOREM 4.2. *Let $0 < c \leq t/\mathbb{E}Z(t)_+$. Then*

$$\mathbb{P}\{Z(t) > c\mathbb{E}Z(t)_+\} \lesssim \frac{t^2}{(\mathbb{E}Z(t)_+)^3}\left(\frac{1}{c^3} + \frac{1}{c^4}\right).$$

PROOF. Note that, for any $\lambda \geq 1$,

$$\mathbb{P}\{Z(t) > u\} = \mathbb{P}\{\exists z > u : N(z) + A_t(z) = L(t,t)\}$$
$$\leq \mathbb{P}\{\exists z > u : N(z) + L_\lambda(t,t) - L_\lambda(z,0) \geq L(t,t)\}$$
$$= \mathbb{P}\{\exists z > u : N(z) - L_\lambda(z,0) \geq L(t,t) - L_\lambda(t,t)\}.$$

Since $L_\lambda(\cdot, 0)$ is a thickening of $L(\cdot, 0)$, we get that $\tilde{N}_{\lambda-1}(z) := L_\lambda(z,0) - N(z)$ is in itself a Poisson process with intensity $\lambda - 1$. This means that

$$\mathbb{P}\{Z(t) > u\} \leq \mathbb{P}\{\tilde{N}_{\lambda-1}(u) \leq L_\lambda(t,t) - L(t,t)\}.$$

To have a useful bound for all $0 \leq u \leq \frac{3}{4}t$, we choose $\lambda$ such that

$$\mathbb{E}\tilde{N}_{\lambda-1}(u) - \mathbb{E}\{L_\lambda(t,t) - L(t,t)\} = (\lambda - 1)u - t\left(\lambda + \frac{1}{\lambda} - 2\right)$$

is maximal. This means that we choose

$$\lambda_u = (1 - u/t)^{-1/2}.$$

Some useful elementary inequalities, that hold for all $0 < u \leq \frac{3}{4}t$, are

$$\lambda_u \leq 2,$$
(4.1) $$\mathbb{E}\tilde{N}_{\lambda_u-1}(u) - \mathbb{E}\{L_{\lambda_u}(t,t) - L(t,t)\} \geq \tfrac{1}{4}u^2/t,$$
$$\lambda_u - 1/\lambda_u \leq 2u/t.$$

Note that, due to (2.7) and (3.6),

$$\text{Var}\{L_{\lambda_u}(t,t) - L(t,t)\} \leq 2(\text{Var}\{L_{\lambda_u}(t,t)\} + \text{Var}\{L(t,t)\})$$
$$\leq 8\mathbb{E}Z(t)_+ + 2t(\lambda_u - 1/\lambda_u)$$
$$\leq 8\mathbb{E}Z(t)_+ + 4u.$$

Now we can use Chebyshev's inequality:

$$\mathbb{P}\{Z(t) > u\} \leq \mathbb{P}\{\tilde{N}_{\lambda_u-1}(u) \leq L_{\lambda_u}(t,t) - L(t,t)\}$$
$$\leq \mathbb{P}\{\tilde{N}_{\lambda_u-1}(u) \leq \mathbb{E}\tilde{N}_{\lambda_u-1}(u) - u^2/(8t)\}$$
$$+ \mathbb{P}\{L_{\lambda_u}(t,t) - L(t,t) \geq \mathbb{E}\tilde{N}_{\lambda_u-1}(u) - u^2/(8t)\}$$

$$
\begin{aligned}
(4.2) \quad &\leq \frac{64t^2(\lambda_u - 1)u}{u^4} \\
&\quad + \mathbb{P}\{L_{\lambda_u}(t,t) - L(t,t) \geq \mathbb{E}\{L_{\lambda_u}(t,t) - L(t,t)\} + u^2/(8t)\} \\
&\lesssim \frac{t^2}{u^3} + \frac{64t^2(8\mathbb{E}Z(t)_+ + 4u)}{u^4} \\
&\lesssim \frac{t^2}{u^3} + \frac{t^2 \mathbb{E}Z(t)_+}{u^4}.
\end{aligned}
$$

If $t \geq u \geq \frac{3}{4}t$, we see that

$$\mathbb{P}\{Z(t) > u\} \leq \mathbb{P}\left\{Z(t) > \frac{3}{4}t\right\} \lesssim \frac{t^2}{u^3} + \frac{t^2 \mathbb{E}Z(t)_+}{u^4},$$

where we use (4.2). This means that (4.2) is true for all $0 \leq u \leq t$. The theorem now follows from choosing $u = c\mathbb{E}Z(t)_+$. □

With this theorem we can show that $\mathbb{E}Z(t)_+ = O(t^{2/3})$.

COROLLARY 4.3.  *Let $Z(t)_+$ and $L(t,t)$ be defined as in (3.6). Then*

$$\limsup_{t \to \infty} \frac{\mathbb{E}Z(t)_+}{t^{2/3}} = \limsup_{t \to \infty} \frac{\mathrm{Var}(L(t,t))}{2t^{2/3}} < +\infty.$$

PROOF. Using (3.6), we only have to prove the statement for $\mathbb{E}Z(t)_+$. Suppose there exists a sequence $t_n \uparrow +\infty$ such that

$$\lim_{n \to \infty} \frac{\mathbb{E}Z(t_n)_+}{t_n^{2/3}} = +\infty.$$

Using Theorem 4.2, we see that

$$\mathbb{P}\{Z(t_n)_+ > c\mathbb{E}Z(t_n)_+\} \lesssim \frac{t_n^2}{(\mathbb{E}Z(t_n)_+)^3}\left(\frac{1}{c^3} + \frac{1}{c^4}\right) \wedge 1.$$

Using dominated convergence [note that $t_n^2/(\mathbb{E}Z(t_n)_+)^3$ is a bounded sequence], this shows that

$$\int_0^\infty \mathbb{P}\{Z(t_n)_+ > c\mathbb{E}Z(t_n)_+\}\,dc \to 0,$$

which would imply the absurd assertion that

$$\mathbb{E}\left\{\frac{Z(t_n)_+}{\mathbb{E}Z(t_n)_+}\right\} \to 0. \qquad \square$$

As a corollary we get the following:

COROLLARY 4.4. *Let $c \geq 1$. Then*
$$\mathbb{P}\{Z(t) > ct^{2/3}\} \lesssim \frac{1}{c^3}.$$

PROOF. This is an immediate consequence of Theorem 4.2 and the previous corollary. □

REMARK. This result can be compared to a result on transversal fluctuations of a longest NE path in [7]. He shows that all longest *stricly* NE paths from $(0,0)$ to $(t,t)$ remain in a strip along the diagonal of width $t^\gamma$, with probability tending to 1, as $t \to \infty$, for any $\gamma > 2/3$. To this end, he uses the analytic results in [2]. Section 3, in combination with Corollary 4.4, shows that the transversal fluctuation at time $s < t$ of a longest *weakly* NE path from $(0,0)$ to $(t,t)$ is of order $(t-s)^{2/3}$. This is due to the fact that any longest weakly NE path lies within the reflected trajectories of a second class particle and a dual second class particle; see Figure 4. Our result on weakly NE paths implies the same result for strictly NE paths, as the following short argument will show: consider the longest strictly NE path that is to the right of all other longest strictly NE paths and suppose, at time $s < t$, this path is to the right of the diagonal. Since the sources and sinks are independent of this event, we have that, given this event, there is still a probability of at least $1/2$ that $Z(t) > 0$. The corresponding longest weakly NE path (so the weakly NE path that is most to the right) cannot be to the left of the considered longest strictly NE path (because they cannot cross), so the order of the fluctuations to the right of a longest strictly NE path cannot be higher than the same order for longest weakly NE paths. For fluctuations to the left, a similar argument holds. The remark at the end of Section 6 discusses the corresponding lower bound result.

**5. Strictly NE paths and restricted weakly NE paths.** To get a lower bound on $\mathbb{E}Z(t)_+$, we need to control the difference between a strictly NE path and a weakly NE path in a stationary Hammersley process $L_1$ (with source intensity 1), where the weakly NE path is only allowed to pick up sources in an interval $[0, \varepsilon t^{2/3}] \times \{0\}$. To do this, we consider another independent Hammersley process $L_\lambda$ on $[0,t]^2$ with source intensity $\lambda$, sink intensity $1/\lambda$ and $\alpha$-intensity 1; for this process, the sources, sinks and $\alpha$-points are independent of the corresponding processes for $L_1$. Coupled to this process $L_\lambda$, we consider $L_0$ as the corresponding (nonstationary) Hammersley process that uses the same $\alpha$-points, but has no sources or sinks.

We denote $L_0(x,t)$ as the number of particles (i.e., the number of crossings of space–time paths) of the Hammersley process without sources or sinks with the segment $[0,x] \times \{t\}$. Note that, for $0 < z < t$,

(5.1) $$A_t(0) - A_t(z) \stackrel{\mathcal{D}}{=} L_0(t,t) - L_0(t-z,t).$$

This follows from the fact that $L_0(t - z, t)$ equals the length of the longest (strictly) NE path from $(0,0)$ to $(t - z, t)$. Also note that $\{N(z) : z \in (0, t)\}$ [the number of sources of the process $L_1$ in the interval $(0, t)$] and $\{L_\lambda(t, t) - L_\lambda(t - z, t) : z \in (0, t)\}$ are two independent Poisson processes.

Define $X'_\lambda(t)$ as the position at time $t$ of a dual second class particle of the process $L_\lambda$, that started in $(0, 0)$. Then we know that

$$(5.2) \quad x < y < X'_\lambda(t) \quad \Longrightarrow \quad L_\lambda(y, t) - L_\lambda(x, t) \leq L_0(y, t) - L_0(x, t).$$

This is due to the fact that if we leave out all the sources of $L_\lambda$, the space–time paths do not change above the trajectory of $X'_\lambda$ (this is one of the key ideas in [4]; the reader might want to check this fact by looking at Figure 2). This means that if we define the process $\underline{L}$ as a Hammersley process that uses the same $\alpha$-points and sinks as $L_\lambda$, but starts without any sources, we have that

$$x < X'_\lambda(t) \quad \Longrightarrow \quad L_\lambda(x, t) = \underline{L}(x, t).$$

Inequality (5.2) now follows from the fact that the set of particles of $\underline{L}$ is at all times a subset of the set of particles of the Hammersley process $L_0$, since this process has no sinks, whereas $\underline{L}$ does have sinks.

THEOREM 5.1. *Fix $L > 0$. Then*

$$\limsup_{t \to \infty} \mathbb{P}\left\{\sup_{z \in [0, \varepsilon t^{2/3}]} \{N(z) + A_t(z)\} - A_t(0) \geq L t^{1/3}\right\} = O(\varepsilon^3), \quad \varepsilon \downarrow 0.$$

PROOF. We will use the auxiliary Hammersley process constructed above. Define

$$\lambda = 1 - r t^{-1/3}.$$

If $X'_\lambda(t) \geq t$, (5.2) tells us that, for $0 < z < t$,

$$L_\lambda(t, t) - L_\lambda(t - z, t) \leq L_0(t, t) - L_0(t - z, t).$$

Using (5.1) and defining

$$\tilde{N}_\lambda(z) = L_\lambda(t, t) - L_\lambda(t - z, t),$$

we see that (remember that $N$ and $A_t$ are independent)

$$(5.3) \quad \mathbb{P}\left\{\sup_{z \in [0, \varepsilon t^{2/3}]} \{N(z) + A_t(z)\} - A_t(0) \geq L t^{1/3}\right\}$$
$$\leq \mathbb{P}\left\{\sup_{z \in [0, \varepsilon t^{2/3}]} \{N(z) - \tilde{N}_\lambda(z)\} \geq L t^{1/3}\right\} + \mathbb{P}\{X'_\lambda(t) < t\}.$$

The second term on the right-hand side of (5.3) can be bounded using (2.5), (3.5) and Corollary 4.4:

$$
\begin{aligned}
\mathbb{P}\{X'_\lambda(t) < t\} &= \mathbb{P}\{X'(t/\lambda) < \lambda t\} \\
&= \mathbb{P}\{Z'(t/\lambda) > t(1/\lambda - \lambda)\} \\
&\leq \mathbb{P}\{Z(t/\lambda) > t(1/\lambda - \lambda)\} \\
&= \mathbb{P}\{Z(t/(1-rt^{-1/3})) > rt^{2/3}(2-rt^{-1/3})/(1-rt^{-1/3})\} \\
&\leq \mathbb{P}\{Z(\tilde{t}) > r\tilde{t}^{\,2/3}\} \lesssim r^{-3},
\end{aligned}
\tag{5.4}
$$

for all $r \in [1, t^{1/3})$, applying Corollary 4.4 with argument $\tilde{t} = t/(1-rt^{-1/3})$.

The first term on the right-hand side of (5.3) concerns a hitting time for the difference of two independent Poisson processes. After rescaling, this can be written as

$$\mathbb{P}\{\exists_{0 \leq z \leq \varepsilon} : t^{-1/3}\{N(zt^{2/3}) - \tilde{N}_\lambda(zt^{2/3})\} \geq L\}.$$

The process $z \mapsto t^{-1/3}\{N(zt^{2/3}) - \tilde{N}_\lambda(zt^{2/3})\}$ converges, as $t \to \infty$, to the drifting Brownian motion process

$$W_r(z) \stackrel{\text{def}}{=} W(2z) + rz, \qquad z \geq 0, \tag{5.5}$$

in the topology of uniform convergence on compacta, where $W$ is standard Brownian motion on $\mathbb{R}_+$. Hence, we get, by a standard application of Donsker's theorem,

$$
\begin{aligned}
\lim_{t \to \infty} \mathbb{P}\{\exists_{0 \leq z \leq \varepsilon} &: t^{-1/3}\{N(zt^{2/3}) - \tilde{N}_\lambda(zt^{2/3})\} \geq L\} \\
&= \mathbb{P}\Big\{\sup_{z \in [0,\varepsilon]} W_r(z) \geq L\Big\}.
\end{aligned}
\tag{5.6}
$$

We now get, for $r < L/\varepsilon$,

$$
\begin{aligned}
\mathbb{P}\Big\{\sup_{z \in [0,\varepsilon]} W_r(z) \geq L\Big\} &\leq \mathbb{P}\Big\{\sup_{z \in [0,2\varepsilon]} W(z) \geq L - \varepsilon r\Big\} \\
&= \mathbb{P}\Big\{\sup_{z \in [0,1]} W(2\varepsilon z)/\sqrt{2\varepsilon} \geq (L - \varepsilon r)/\sqrt{2\varepsilon}\Big\} \\
&= \mathbb{P}\Big\{\sup_{z \in [0,1]} W(z) \geq (L - \varepsilon r)/\sqrt{2\varepsilon}\Big\} \\
&= \sqrt{\frac{2}{\pi}} \int_{(L-\varepsilon r)/\sqrt{2\varepsilon}}^\infty e^{-(1/2)u^2} \, du.
\end{aligned}
\tag{5.7}
$$

Taking $r = L/(2\varepsilon)$, we get

$$\sqrt{\frac{2}{\pi}} \int_{(L-\varepsilon r)/\sqrt{2\varepsilon}}^{\infty} e^{-(1/2)u^2} \, du$$

$$= 2 \int_{L/\sqrt{8\varepsilon}}^{\infty} \frac{e^{-(1/2)u^2}}{\sqrt{2\pi}} \, du \sim \frac{2\sqrt{8\varepsilon}\, e^{-L^2/(16\varepsilon)}}{L\sqrt{2\pi}}, \qquad \varepsilon \downarrow 0,$$

using Mills' ratio approximation for the tail of a normal distribution in the last step. This means that, with this choice of $r$, our estimate for the second term on the right-hand side of (5.3) is dominant, so (5.4) now proves the theorem. $\square$

Note that we could use the proof of the theorem to show that $L$ can even go to 0 at a certain speed when $\varepsilon \to 0$, and still the considered probability would go to 0, uniformly in $t$.

**6. Lower bound for $\mathbb{E}Z(t)_+$.** We wish to bound the probability that $Z(t) \in [0, \varepsilon t^{2/3}]$. In order to do this, we again introduce an independent auxiliary stationary Hammersley process $L_\lambda$, but now with intensity $\lambda > 1$. In fact, we will choose

$$\lambda = 1 + rt^{-1/3}.$$

Coupled to this process, we again consider the Hammersley process $L_0$ without sources or sinks, but with the same $\alpha$-points. This time, however, we will leave out the sinks of the stationary process; to be more precise, we have that

(6.1) $\quad y \geq x > X_\lambda(t) \quad \implies \quad L_0(y,t) - L_0(x,t) \leq L_\lambda(y,t) - L_\lambda(x,t).$

The reason is that if we consider the stationary process $L_\lambda$ and leave out the sinks of this process, below the trajectory of $X_\lambda$ the space–time paths do not change. The inequality then follows from the fact that the set of particles of a process which, like $L_0$, has no sinks, but starts with sources and uses the same $\alpha$-points, will at all times be a superset of the set of particles of the Hammersley process $L_0$. Compare this to the explanation of (5.2).

We again define

$$\tilde{N}_\lambda(z) = L_\lambda(t,t) - L_\lambda(t-z,t)$$

and will show the following result.

THEOREM 6.1.

$$\lim_{\varepsilon \downarrow 0} \limsup_{t \to \infty} \mathbb{P}\{0 \leq Z(t) \leq \varepsilon t^{2/3}\} = 0.$$

PROOF. Let $\eta > 0$. It is enough to find $\varepsilon > 0$ such that
$$\limsup_{t\to\infty} \mathbb{P}\{0 \leq Z(t) \leq \varepsilon t^{2/3}\} < 3\eta.$$

For any $L, r > 0$ and $\lambda = 1 + rt^{-1/3}$, we have

$$\mathbb{P}\{0 \leq Z(t) \leq \varepsilon t^{2/3}\} \leq \mathbb{P}\left\{\sup_{z > \varepsilon t^{2/3}}(N(z) + A_t(z)) < \sup_{z \in [0, \varepsilon t^{2/3}]}(N(z) + A_t(z))\right\}$$

$$(6.2) \qquad \leq \mathbb{P}\left\{\sup_{z > \varepsilon t^{2/3}}\{N(z) - (A_t(0) - A_t(z))\} < Lt^{1/3}\right\}$$

$$+ \mathbb{P}\left\{\sup_{z \in [0, \varepsilon t^{2/3}]}\{N(z) - (A_t(0) - A_t(z))\} > Lt^{1/3}\right\},$$

since for any $L \in \mathbb{R}$, $X < Y$ implies that either $X < L$ or $Y > L$, so
$$\mathbb{P}(X < Y) \leq \mathbb{P}(\{X < L\} \cup \{Y > L\}) \leq \mathbb{P}(X < L) + \mathbb{P}(L < Y),$$

which allows us to "optimize" over $L$. For the first term of (6.2), we want to use (6.1), which implies that, on the event $\{X_\lambda(t) \leq t - rt^{2/3}\}$, we know that, for all $z \leq rt^{2/3}$,
$$L_0(t, t) - L_0(t - z, t) \leq L_\lambda(t, t) - L_\lambda(t - z, t).$$

Therefore,

$$\mathbb{P}\{0 \leq Z(t) \leq \varepsilon t^{2/3}\}$$

$$\leq \mathbb{P}\left\{\sup_{z \in [\varepsilon t^{2/3}, rt^{2/3}]}\{N(z) - \tilde{N}_\lambda(z)\} < Lt^{1/3}\right\} + \mathbb{P}\{X_\lambda(t) > t - rt^{2/3}\}$$

$$(6.3) \qquad + \mathbb{P}\left\{\sup_{z \in [0, \varepsilon t^{2/3}]}\{N(z) - (A_t(0) - A_t(z))\} > Lt^{1/3}\right\}$$

$$= \mathbb{P}\left\{\sup_{z \in [\varepsilon t^{2/3}, rt^{2/3}]}\{N(z) - \tilde{N}_\lambda(z)\} < Lt^{1/3}\right\} + \mathbb{P}\{X(t/\lambda) > \lambda t - \lambda r t^{2/3}\}$$

$$+ \mathbb{P}\left\{\sup_{z \in [0, \varepsilon t^{2/3}]}\{N(z) - (A_t(0) - A_t(z))\} > Lt^{1/3}\right\}.$$

Using (3.5) (note that $\lambda t - \lambda r t^{2/3} \geq t/\lambda$ when $r \leq \frac{1}{2}t^{1/3}$), $Z'(t) \leq Z(t)$ and $-Z'(t) \stackrel{\mathcal{D}}{=} Z(t)$, we get, for the second term in (6.3),

$$\mathbb{P}\{X(t/\lambda) > \lambda t - \lambda r t^{2/3}\} = \mathbb{P}\{Z(\lambda t - \lambda r t^{2/3}) < (1/\lambda - \lambda)t + \lambda r t^{2/3}\}$$

$$\leq \mathbb{P}\{Z'(\lambda t - \lambda r t^{2/3}) < (1/\lambda - \lambda)t + \lambda r t^{2/3}\}$$

$$= \mathbb{P}\{-Z'(\lambda t - \lambda r t^{2/3}) > (\lambda - 1/\lambda)t - \lambda r t^{2/3}\}$$

$$= \mathbb{P}\{Z(\lambda t - \lambda r t^{2/3}) > (\lambda - 1/\lambda)t - \lambda r t^{2/3}\}$$

$$= \mathbb{P}\left\{Z(t(1 - r^2 t^{-2/3})) > \frac{r t^{2/3}}{1 + r t^{-1/3}} - r^2 t^{1/3}\right\}.$$

This means that we can start by choosing $r$ sufficiently large to ensure that the second term is smaller than $\eta$, since

$$\mathbb{P}\left\{Z(t(1 - r^2 t^{-2/3})) > \frac{r t^{2/3}}{1 + r t^{-1/3}} - r^2 t^{1/3}\right\}$$

$$= \mathbb{P}\left\{Z(t + O(t^{1/3})) > r t^{2/3} + O(t^{1/3})\right\} = O(r^{-3}), \qquad t \to \infty,$$

where we use Corollary 4.4 in the last step.

Now we turn to the first term in (6.3). This term is very similar to the term we found in the previous section. We have, as in the proof of Theorem 5.1, that the process

$$z \mapsto t^{-1/3}\{N(zt^{2/3}) - \tilde{N}_\lambda(zt^{2/3})\}$$

converges, as $t \to \infty$, to a drifting Brownian motion process

$$(6.4) \qquad W_r(z) \stackrel{\text{def}}{=} W(2z) - rz, \ z \geq 0,$$

in the topology of uniform convergence on compacta, where $W$ is standard Brownian motion on $\mathbb{R}_+$ (this time the drift is negative instead of positive). Hence, we get, again using Donsker's theorem,

$$\lim_{t \to \infty} \mathbb{P}\left\{\sup_{z \in [\varepsilon t^{2/3}, r t^{2/3}]} \{N(z) - \tilde{N}_\lambda(z)\} \leq L t^{1/3}\right\}$$

$$= \mathbb{P}\left\{\sup_{z \in [\varepsilon, r]} W_r(z) \leq L\right\}$$

(6.5)

$$= \mathbb{P}\left\{\sup_{z \in [0,r]} W_r(z) \leq L\right\} + \mathbb{P}\left\{\sup_{z \in [\varepsilon, r]} W_r(z) \leq L, \sup_{z \in [0,\varepsilon]} W_r(z) > L\right\}$$

$$\leq \mathbb{P}\left\{\sup_{z \in [0,r]} W_r(z) \leq L\right\} + \mathbb{P}\left\{\sup_{z \in [0,\varepsilon]} W_r(z) > L\right\}.$$

Since

$$\lim_{L \downarrow 0} \mathbb{P}\left\{\sup_{z \in [0,r]} W_r(z) \leq L\right\} = 0,$$

we can choose $L = L(\eta) > 0$ sufficiently small to ensure

$$\mathbb{P}\left\{\sup_{z \in [0,r]} W_r(z) \leq L\right\} < \eta/2.$$

It is also clear from the argument of the proof of Theorem 5.1 [see (5.7)] that we can next choose $\varepsilon > 0$ sufficiently small to ensure that

$$\mathbb{P}\left\{\sup_{z \in [0,\varepsilon]} W_r(z) > L\right\} \leq \mathbb{P}\left\{\sup_{z \in [0,\varepsilon]} W(2z) > L\right\} < \eta/2,$$

for this choice of $L = L(\eta) > 0$.

It is now seen from (6.5) that this bounds the first term of (6.3) from above by $\eta$ [remember that we have already fixed $r > 0$ to bound the second term of (6.3)]. Finally, we can choose $\varepsilon > 0$ so small that the third term in (6.3) is smaller than $\eta$, using Theorem 5.1. This completes the proof. $\square$

COROLLARY 6.2. *Let $Z(t)_+$ and $L(t,t)$ be defined as in* (3.6). *Then*

$$\liminf_{t \to \infty} \frac{\mathbb{E}Z(t)_+}{t^{2/3}} = \liminf_{t \to \infty} \frac{\mathrm{Var}(L(t,t))}{2t^{2/3}} > 0.$$

PROOF. Using (3.6), we only have to prove the statement for $\mathbb{E}Z(t)_+$. Suppose $t_n \to \infty$ such that

$$\frac{\mathbb{E}Z(t_n)_+}{t_n^{2/3}} \to 0.$$

Then

$$\mathbb{P}\{Z(t_n) > \varepsilon t_n^{2/3}\} \leq \frac{\mathbb{E}Z(t_n)_+}{\varepsilon t_n^{2/3}} \to 0.$$

Since $-Z'(t) \stackrel{\mathcal{D}}{=} Z(t)$ and $Z'(t) \leq Z(t)$, we have that $\mathbb{P}\{Z(t) \geq 0\} \geq 1/2$, for all $t > 0$. This would mean that, for all $\varepsilon > 0$,

$$\liminf_{n \to \infty} \mathbb{P}\{0 \leq Z(t_n) \leq \varepsilon t_n^{2/3}\} \geq \tfrac{1}{2},$$

which would contradict Theorem 6.1. $\square$

REMARK. We can again make a comparison with the results in [7]. Johansson shows that the probability that all longest *strictly* NE paths from $(0,0)$ to $(t,t)$ stay within a strip around the diagonal of width $t^\gamma$ *does not* tend to 1, as $t \to \infty$, for all $\gamma < 2/3$, again using the analytical results in [2]. Our Theorem 6.1 shows that the probability that all longest *weakly* NE paths from $(0,0)$ to $(t,t)$ stay within a strip around the diagonal of width $t^\gamma$ *tends to* 0, as $t \to \infty$, for all $\gamma < 2/3$ (see also the Remark at the end of Section 4).

**7. Tightness results.** In the preceding sections it was shown, using the "hydrodynamical methods" of [4] that, for a stationary version of Hammersley's process, with intensity 1 for the Poisson point processes on the axes and in the plane, the variance of the length of a longest weakly NE path $L(t,t)$ is of order $t^{2/3}$, in the sense that

$$(7.1) \quad 0 < \liminf_{t \to \infty} t^{-2/3} \operatorname{Var}(L(t,t)) \leq \limsup_{t \to \infty} t^{-2/3} \operatorname{Var}(L(t,t)) < \infty.$$

This means, in particular, that, for any $t > 0$, the sequence

$$n^{-1/3}\{L(nt, nt) - 2nt\}, \qquad n = 1, \ldots,$$

is tight.

As noted in [9], the distributional limit result for $n^{-1/3}\{L_0(nx, nt) - 2n\sqrt{xt}\}$ for Hammersley's process without sources and sinks in [2] can be translated into a limit of $Y_n/n^{1/3}$, where

$$(7.2) \qquad Y_n \stackrel{\text{def}}{=} z_{nt}([nx]) - nx^2/(4t),$$

and $z_{nt}([nx])$ is the $[nx]$th particle at time $nt$, counting particles at time $nt$ from the left. Theorem 3.2 in [9] gives a tightness result for a more general version of $Y_n$, in the context of a version of Hammersley's process on the whole line, with a (possibly) random initial state. The result is that, under his conditions D and E, the sequence

$$Y_n/(n^{1/3} \log n), \qquad n = 1, \ldots$$

is tight. He conjectures that, in fact, $n^{1/3} \log n$ can be replaced by $n^{1/3}$. The results, derived above, are a further indication that indeed $n^{1/3} \log n$ might be replaced by $n^{1/3}$, and that this can be derived by hydrodynamical methods.

For the stationary version of Hammersley's process, with intensities 1 of the Poisson processes in the plane and on the axes, we can define $z_{nt}([2nt])$ as the location of the $[2nt]$th source at time $nt$, where we count the sources from left to right, starting with the first source to the right of zero. Note that at time zero the particles are just the sources. The particles, escaping through a sink, are given location zero at times larger than or equal to the time of escape.

With this definition, our results give tightness of the sequences

$$(7.3) \qquad n^{-1/3}\{z_{nt}([2nt]) - nt\}, \qquad n = 1, 2, \ldots,$$

for each $t > 0$. This can be seen in the following way. We have, for $M > 0$ and $t > 0$, the "switch relation"

$$(7.4) \quad n^{-1/3}\{z_{nt}([2nt]) - nt\} > M \quad \Longleftrightarrow \quad L(nt + Mn^{1/3}, nt) < [2nt].$$

Theorem 2.1 yields

$$\begin{aligned}\operatorname{Var}(L(nt+Mn^{1/3},nt)) &= -nt-Mn^{1/3}+nt+2\mathbb{E}(nt+Mn^{1/3}-X_1(nt))_+\\ &= -Mn^{1/3}+2\mathbb{E}(nt+Mn^{1/3}-X_1(nt))_+\\ &\leq Mn^{1/3}+2\mathbb{E}(nt-X_1(nt))_+,\end{aligned}$$

where we use $(Y+Z)_+ \leq Y_+ + Z_+$ in the last step. Theorem 2.1, applied in the opposite direction, yields

$$2\mathbb{E}(nt-X_1(nt))_+ = \operatorname{Var}(L(nt,nt)).$$

Hence, we get, by (7.1) and Chebyshev's inequality,

$$\begin{aligned}&\mathbb{P}\{n^{-1/3}\{z_{nt}([2nt])-nt\} > M\}\\ &= \mathbb{P}\{L(nt+Mn^{1/3},nt)-2nt-Mn^{1/3} < [2nt]-2nt-Mn^{1/3}\}\\ &\leq \frac{Mn^{1/3}+\operatorname{Var}(L_1(nt,nt))}{\{Mn^{1/3}+2nt-[2nt]\}^2} \lesssim \frac{Mn^{1/3}+O((nt)^{2/3})}{M^2 n^{2/3}} = O(M^{-1}).\end{aligned}$$

We similarly get

$$\mathbb{P}\{n^{-1/3}\{z_{nt}([2nt])-nt\} \leq -M\} = O(M^{-1}),$$

using that, if $nt - Mn^{-1/3} > 0$,

(7.5) $\quad n^{-1/3}\{z_{nt}([2nt])-nt\} \leq -M \iff L(nt-Mn^{1/3},nt) \geq [2nt],$

which proves the tightness of the sequences (7.3).

Although the tightness of the sequence $(Y_n/n^{1/3})$ for the Hammersley process without sources or sinks, as defined in (7.2), is known from the results of [2], it is of some interest to derive this from the results of the preceding sections. The tightness will follow from

(7.6) $\quad n^{-1/3}\{L_0(nx,nt) - 2\sqrt{xt}\} = O_p(1), \qquad n \to \infty,$

for all $x, t > 0$, where $L_0(nx, xt)$ is a strictly NE path from $(0,0)$ to $(nx, xt)$, and where the intensity of the Poisson process in the first quadrant is equal to 1.

We again have a "switch relation" similar to (7.4):

$$n^{-1/3}\{z_{nt}([2nt])-nt\} > M \iff L_0(nt+Mn^{1/3},nt) < [2nt].$$

From [1], we know that

$$L_0(nx,nt) \stackrel{\mathcal{D}}{=} L_0(n\sqrt{xt}, n\sqrt{xt}),$$

so if we show that, for each $t > 0$,

(7.7) $\quad n^{-1/3}\{L_0(nt,nt) - 2nt\} = O_p(1), \qquad n \to \infty,$

we get, for each $\varepsilon > 0$ and $t > 0$,

$$\mathbb{P}\{n^{-1/3}\{z_{nt}([2nt]) - nt\} > M\}$$
$$= \mathbb{P}\{L_0(nt + Mn^{1/3}, nt) - 2nt - Mn^{1/3} < [2nt] - 2nt - Mn^{1/3}\}$$
$$= \mathbb{P}\{L_0(nt\sqrt{1 + Mn^{-2/3}/t}, nt\sqrt{1 + Mn^{-2/3}/t})$$
$$- 2nt\sqrt{1 + Mn^{-2/3}/t} + O(M^2 n^{-1/3})$$
$$< [2nt] - 2nt - Mn^{1/3}\}$$
$$< \varepsilon,$$

for sufficiently large $M = M(\varepsilon) > 0$ and all $n \geq n_0(M, \varepsilon)$. Relation (7.7) similarly implies

$$\mathbb{P}\{n^{-1/3}\{z_{nt}([2nt]) - nt\} \leq -M\} < \varepsilon,$$

for sufficiently large $M = M(\varepsilon) > 0$ and all $n \geq n_0(M, \varepsilon)$, using

$$n^{-1/3}\{z_{nt}([2nt]) - nt\} \leq -M \iff L_0(nt - Mn^{1/3}, nt) \geq [2nt].$$

In order to prove (7.7), it is sufficient to show

(7.8) $$t^{-1/3}\{L_0(t, t) - 2t\} = O_p(1), \qquad t \to \infty.$$

Now first note that the length $L_0(t,t)$ of a longest strictly NE path from $(0,0)$ to $(t,t)$ is the same as $A_t(0)$ in the proof of Theorem 5.1. Let, for $\lambda = 1 - rt^{-1/3}$, $L_\lambda$ be defined as in the proof of Theorem 5.1. By (5.3), we have

(7.9)
$$\mathbb{P}\left\{\sup_{z \in [0, Kt^{2/3}]} \{N(z) + A_t(z)\} - A_t(0) \geq Lt^{1/3}\right\}$$
$$\leq \mathbb{P}\left\{\sup_{z \in [0, Kt^{2/3}]} \{N(z) - \tilde{N}_\lambda(z)\} \geq Lt^{1/3}\right\} + \mathbb{P}\{X'_\lambda(t) < t\}.$$

We first deal with the second term on the right-hand side of (7.9). By (5.4), we have

(7.10) $$\mathbb{P}\{X'_\lambda(t) < t\} = O(r^{-3}),$$

uniformly for $r \in [1, \tfrac{1}{2}t^{1/3}]$. To deal with the first term on the right-hand side of (7.9), we first note that

$$z \mapsto M_r(z) \stackrel{\text{def}}{=} t^{-1/3}\{N(zt^{2/3}) - \tilde{N}_\lambda(zt^{2/3})\} - rz, \qquad z \geq 0,$$

is a zero-mean martingale. So we get

$$\mathbb{P}\left\{\sup_{z\in[0,Kt^{2/3}]}\{N(z)-\tilde{N}_\lambda(z)\}\geq Lt^{1/3}\right\}=\mathbb{P}\left\{\sup_{z\in[0,K]}\{M_r(z)+rz\}\geq L\right\}$$

$$\leq \mathbb{P}\left\{\sup_{z\in[0,K]}M_r(z)\geq L-rK\right\}.$$

Taking $r=L/(2K)$ and using Doob's submartingale inequality, we get

$$\mathbb{P}\left\{\sup_{z\in[0,Kt^{2/3}]}\{N(z)-\tilde{N}_\lambda(z)\}\geq Lt^{1/3}\right\}\leq \mathbb{P}\left\{\sup_{z\in[0,K]}M_r(z)\geq L/2\right\}$$

$$\leq \mathbb{P}\left\{\sup_{z\in[0,K]}M_r(z)^2\geq L^2/4\right\}$$

$$\leq \frac{4EM_r(K)^2}{L^2}\lesssim K/L^2.$$

We also have, by Corollary 4.4,

$$\mathbb{P}\left\{\sup_{z\in[-Kt^{2/3},Kt^{2/3}]}\{N(z)+A_t(z)\}\neq \sup_{z\in[-t,t]}\{N(z)+A_t(z)\}\right\}$$

$$\leq 2\mathbb{P}\{Z(t)>Kt^{2/3}\}\lesssim 1/K^3.$$

So, taking $K=L^{7/12}$ [note that this means that $r=L/(2K)=\frac{1}{2}L^{5/12}\leq \frac{1}{2}t^{1/3}$, for $L\leq t^{4/5}$], we obtain

$$\mathbb{P}\left\{\sup_{z\in[-t,t]}\{N(z)+A_t(z)\}-A_t(0)\geq Lt^{1/3}\right\}$$

$$=\mathbb{P}\{L_1(t,t)-A_t(0)\geq Lt^{1/3}\}\lesssim L^{-5/4},$$

for all $L\leq t^{4/5}$.

If $L>t^{4/5}$, we first note that

$$\mathbb{P}\{L_1(t,t)-A_t(0)\geq Lt^{1/3}\}\leq \mathbb{P}\{L_1(t,t)\geq Lt^{1/3}\}$$

$$\leq 2\mathbb{P}\{P_t\geq \tfrac{1}{2}Lt^{1/3}\}\leq 2\mathbb{P}\{P_t\geq \tfrac{1}{2}L^{1/6}t\},$$

where $P_t$ is a Poisson variable with expectation $t$. Let $[x]$ denote the largest integer $\leq x$ and let $a=\frac{1}{2}L^{1/6}$. Then, using the Lagrange remainder term in an expansion of $e^t$, we get, for a $\theta\in(0,1)$,

$$\mathbb{P}\{P_t\geq \tfrac{1}{2}L^{1/6}t\}\leq \mathbb{P}\{P_t\geq [at]\}=\frac{t^{[at]}e^{-(1-\theta)t}}{[at]!}\leq \frac{t^{[at]}}{[at]!}.$$

Stirling's formula for the gamma function $\Gamma(x)$ yields that, uniformly in $t \geq 1$,

$$\frac{t^{at}}{\Gamma(at+1)} \sim \frac{1}{\sqrt{2\pi at}} e^{-a(\log a - 1)t}, \qquad a \to \infty.$$

This implies that $\mathbb{P}\{P_t \geq \frac{1}{2}Lt^{1/3}\}$ tends to zero faster than any negative power of $L$, if $L > t^{4/5}$, uniformly in all large $t$ and hence, we can conclude that

$$\mathbb{P}\{L_1(t,t) - A_t(0) \geq Lt^{1/3}\} = O(L^{-5/4}),$$

for all $L \geq 1$, implying

$$\begin{aligned}
0 \leq 2t - \mathbb{E}A_t(0) &= \mathbb{E}\{L_1(t,t) - A_t(0)\} \\
&\lesssim t^{1/3}\left\{1 + \int_1^\infty L^{-5/4}\,dL\right\} = O(t^{1/3}).
\end{aligned} \qquad (7.11)$$

Thus,

$$\mathbb{E}|A_t(0) - 2t| = \mathbb{E}|L_0(t,t) - 2t| = O(t^{1/3}).$$

This proves (7.8) and, as noted above, (7.6) now also follows. This, in turn, proves tightness of the sequence $(Y_n/n^{1/3})$, for $Y_n$ defined by (7.2) for Hammersley's process, starting with the empty configuration on the axes.

Notice that we also proved at the same time

$$\mathbb{E}L_0(nx, nt) = \mathbb{E}L_0(n\sqrt{xt}, n\sqrt{xt}) = \mathbb{E}A_{n\sqrt{xt}}(0) = 2n\sqrt{xt} + O(n^{1/3}),$$

for all $x, t > 0$; see (7.11).

DEPARTMENT OF APPLIED MATHEMATICS (DIAM)
DELFT UNIVERSITY OF TECHNOLOGY
MEKELWEG 4
2628 CD DELFT
THE NETHERLANDS
E-MAIL: e.a.cator@ewi.tudelft.nl
        p.groeneboom@ewi.tudelft.nl